\documentclass[10pt]{amsart}

\usepackage{amsmath,amssymb}
\usepackage{amsfonts}

\theoremstyle{plain} 

\newtheorem{theorem}{Theorem}
\newtheorem*{thm }{Theorem}
\newtheorem{lemma}{Lemma}
\newtheorem{remark}{Remark}
  
\newtheorem{corollary}{Corollary}
\newtheorem*{cor }{Corollary}

\begin{document}
\title{Weak Hyperbolicity on periodic orbits for polynomials.}

\author[J. Rivera-Letelier]{Juan Rivera-Letelier}
\address{J. Rivera-Letelier \\
         Mathematics Department \\
         SUNY at Stony Brook \\
         Stony Brook, NY 11794-3660}
\email{rivera@@math.sunysb.edu}

\date{October 8, 2001}

\newcommand{\diam}{\mbox{\rm diam}}
\newcommand{\dist}{\mbox{\rm dist}}

\def\proof{\par\noindent {\bf Proof.} }

\def\tl{\tilde}
\def\ov{\overline}
\def\cal{\mathcal}

\def\e{\varepsilon}
\def\la{\lambda}

\def\CC{\Bbb C}
\def\DD{\Bbb D}
\def\FF{\Bbb F}
\def\HH{\Bbb H}
\def\RR{\Bbb R}

\begin{abstract}{
	We prove that if the multipliers of the repelling periodic orbits of a complex polynomial grow at least like $n^{5 + \varepsilon}$, for some $\varepsilon > 0$, then the Julia set of the polynomial is locally connected when it is connected.
	As a consequence for a polynomial the presence of a Cremer cycle implies the presence of a sequence of repelling periodic orbits with ``small'' multipliers.
	Somehow surprisingly the proof is based in measure theorical considerations.
}\end{abstract}

\maketitle

\

	Consider a polynomial $P$ with complex coefficients.
	Given a periodic point $p \in \CC$ of $P$ of minimal period $n \ge 1$,
we call
$
	\lambda = (P^n)'(p)
$
the {\it multiplier} of $p$.
	We say that $p$ is {\it repelling}, {\it indifferent} or {\it attracting} if $|\lambda| > 1$, $|\lambda| = 1$ or $|\lambda| < 1$, respectively.


	Recall that the set,
$$
	K(P) = \{ z \in \CC \mid \{ P^n(z) \}_{n \ge 1} \mbox{ is bounded} \},
$$
is called the {\it filled-in Julia set} of $P$ and its boundary is called the {\it Julia set} of $P$, which is denoted by $J(P)$.

	Our main result is the following.

\

\begin{theorem}\label{loc}
	Let $P \in \CC[z]$ be a polynomial with connected Julia set.
	Suppose that there are constants $C > 0$ and $\varepsilon > 0$ such that for every repelling periodic point $p \in \CC$ of period $n$,
$$
	|(P^n)'(p)| \ge C n^{5 + \varepsilon}.
$$
	Then the Julia set of $P$ is locally connected.
\end{theorem}

\


	In fact we prove the stronger statement that $\overline{\CC} - K(P)$ is an integrable domain in the sense of \cite{GS}.

	The class of polynomial (and rational) maps for which the multipliers grow exponentially in the period is studied in \cite{PRS} and the proof of Theorem \ref{loc} is based on a variant of Lemma 2.1 of that paper.

	Recall that a {\it Cremer cycle} is a cycle for which the polynomial is not locally linearizable and whose multiplier is not a root of unity.
	By a theorem of A. Douady and D. Sullivan the Julia set of polynomial with a Cremer cycle is not locally connected, see \cite{Su}.
	So the following corollary follows directly from Theorem \ref{loc}.

\

\begin{corollary}\label{Cremer}
	Let $P \in \CC[z]$ be a polynomial with connected Julia set and with a Cremer cycle.
	Then for every $\varepsilon > 0$ there is a sequence $\{ p_k \}$ of repelling periodic points of $P$ of period $n_k$ satisfying,
$$
	|(P^{n_k})'(p_k)| \le n_k^{5 + \varepsilon}.
$$
\end{corollary}

\

	The corollary remains true if the Julia set is disconnected (see Remark \ref{remark}) and there is an analogous statement for general rational maps, see Remark \ref{remark2}.
	Unfortunately this result does not give information about the location of these periodic points.
	Under certain conditions Cremer cycles imply the existence of the so called ``small'' cycles; see \cite{Y} and \cite{P-M}.

\

\begin{remark}
		It follows by \cite{BvS} (Corollary 1.1) that the hypothesis of Theorem \ref{loc} are satisfied for polynomials such that for every critical value $v \in J(P)$ the derivatives $|(P^n)'(v)|$ grow at least as $n^{\alpha}$, for some $\alpha > 1$ only depending on the degree of $P$.
	Such polynomials also satisfy the so called summability condition, see also \cite{NvS}, \cite{Pr1}, \cite{GS}, \cite{PrUparabolic} and \cite{R-L}.
\end{remark}
\begin{remark}\label{remark}
	To see that Corollary \ref{Cremer} applies for polynomials with disconnected Julia set we first observe that the proof of the theorem applies for polynomial like maps in the sense of \cite{DH}.
	Then if $P \in \CC[z]$ is a polynomial with a Cremer periodic point $p \in \CC$ of period $q \ge 1$, then the restriction of $P^q$ to a suitable neighborhood of $p$ is a polynomial like map with connected Julia set. 
\end{remark}
\begin{remark}\label{remark2}
			A similar method allows to prove that if a rational map $R \in \CC(z)$ has a Cremer cycle, then there is a sequence of periodic points $\{ p_k \}$ of period $n_k$, whose multiplier is bounded by $\exp( C \sqrt{n_k} (\ln n_k)^{\frac{3}{2} + \frac{2}{h}})$, where $C > 0$ is a universal constant and $h > 0$ is any number satisfying $h < HD_{hyp}(R)$.
	The latter is the supremum over the Hausdorff dimensions of hyperbolic sets of $R$.
\end{remark}

\section{Proof of the theorem.}
	Fix a polynomial $P \in \CC[z]$ of degree $d > 1$ and with connected Julia set.
	Consider a base point $w_0 \in \CC - K(P)$ to be chosen in Lemma \ref{expansion} below.

	Theorem \ref{loc} will be reduced to the following lemma, see \cite{GS} for its proof.

\

\begin{lemma}\label{integrable}
	Let $\{ \omega_n \}_{n \ge 1}$ be an increasing sequence such that $\sum_{n \ge 1} \omega_n^{-1} < \infty$.
	If for every $n \ge 1$ and every $z \in P^{-n}(w_0)$ we have $|(P^n)'(z)| \ge \omega_n$, then $J(P)$ is locally connected.
\end{lemma}

\

	Let $\{ \lambda_n \}_{ n \ge 1}$ be an increasing sequence and suppose that for every $n \ge 1$ the repelling periodic points of period $n$ have multipliers of norm at least $\lambda_n > 0$.

	The following lemma estimates the derivative at a good portion of preimages of a base point $w_0 \in \CC - K(P)$.
	Lemmas \ref{distortion0} and \ref{distortion} are distortion lemmas that will allow us to estimate derivatives at {\it all} preimages of $w_0$, as required by Lemma \ref{integrable}.

\

\begin{lemma}\label{expansion}
	Let $\mu$ be an invariant probability measure with positive entropy $h_\mu$ supported on a hyperbolic set $K \subset J(P)$.	
	Then there is a base point $w_0 \in \CC - K(P)$, a constant $C_0 = C_0(w_0) > 0$ and a non-decreasing sequence of integers $\{ \ell(k) \}_{k \ge 1}$ such that the following properties hold.
	\begin{enumerate}
		\item
			There is $L = L(K) \ge 1$ such that $\ell(k + 1) - \ell(k) \le L$, for $k \ge 1$.
			Moreover for every $\varepsilon_0 > 0$ we have $\ell(k) \le (h_\mu^{-1} + \varepsilon_0) \ln k$, for $k$ big enough.
		\item
			For every $k \ge 1$ there is a point $x \in P^{- \ell(k)}(w_0)$ such that for every $w' \in P^{- k}(x)$,
$$
	|(P^{k + \ell(k)})'(w')| \ge C_0 \lambda_{k + \ell(k)}.
$$
	\end{enumerate}
\end{lemma}
{\bf Proof.}
	Since $K$ is a hyperbolic set the following well-known univalent pull-back property holds.
	There is $\delta > 0$ such that for every $z_0 \in K$ there is $\zeta_0 \in K$ such that for every $r > 0$ small there is a univalent pull-back $V \subset B(z_0, r)$ of $B(\zeta_0, 3\delta)$ whose diameter is comparable to $r$.
	Let $R > 1$ such that for any $r > 0$ small and any univalent pull-back of $B(z_0, Rr)$, the corresponding pull-back of $B(z_0, r)$ has diameter at most $\delta$.

	Note that for small $r > 0$ any pull-back of $B(z_0, r)$ will be contained in a definite neighborhood $U$ of $J(P)$.
	Let $m_0 \ge 1$ be such that for every $\zeta_0 \in J(P)$ the set $\cup_{0 \le m \le m_0} P^m(B(\zeta_0, \delta))$ contains $U$.

$1.-$	By Ruelle's inequality $\chi_\mu = \int |P'| d \mu \ge \frac{1}{2} h_\mu > 0$ is positive and then $HD(\mu) = \frac{h_\mu}{\chi_\mu}$ is also positive, where $HD(\mu)$ denotes the infimum of the Hausdorff dimensions of sets $X \subset K$ such that $\mu(X) = 1$; see \cite{Ma} or \cite{PU}.

$2.-$	Choose $\alpha > HD(\mu)^{-1}$.
	By a Borel-Cantelli argument for every $z_0 \in K$, outside an exceptional set of Hausdorff dimension at most $\alpha^{-1} < HD(\mu)$, there is $k_0 = k_0(z_0)$ such that for all $k \ge k_0$ the ball $B(z_0, Rk^{-\alpha})$ is disjoint form $\cup_{0 \le i \le k + m_0} P^i({\rm Crit})$; here ${\rm Crit} \subset \CC$ denotes the set of critical points of $P$.

	So for big values of $k$ all pull-backs of $B(z_0, Rk^{-\alpha})$ by $f^{k + m_0}$ are univalent.

$3.-$	Since $\alpha^{-1} < HD(\mu)$ it follows that the set of such points $z_0 \in K$ has full $\mu$ measure.
	On the other hand it follows by Birkhoff Ergodic Theorem there is a set of $\mu$ positive measure of points in $K$ whose Lyapunov exponent is at least $\chi_\mu = \int |P'| d \mu$, cf. \cite{PU}.

	Thus there is such a point $z_0 \in K$ having Lyapunov exponent at least $\chi_\mu$.

$4.-$	Let $\zeta_0 \in K$ be the corresponding point as explained above and choose any point $w_0 \in B(\zeta_0, \delta) - K(P)$ as a base point.

     	For $k \ge 1$ let $\ell(k)$ be an integer such that there is a univalent pull-back $V \subset B(z_0, k^{-\alpha})$ of $B(\zeta_0, 3\delta)$ by $P^{\ell(k)}$ with $\diam(V)$ comparable to $k^{-\alpha}$.
	We assume that $\{ \ell(k) \}_{k \ge 1}$ is a non-decreasing sequence.
	Note that $\ell(k + 1) - \ell(k) \le L$, for $L \ge 1$ only depending on $K$.

$5.-$ 	Choose $\lambda \in (1, \exp(\chi_\mu))$.
	Since the Lyapunov exponent of $z_0$ is at least $\chi_\mu$ it follows that for $k$ big, $\ell(k) \le \frac{\alpha}{\ln \lambda} \ln k$.
	Moreover for every $\varepsilon_0 > 0$ we can choose $\alpha$ and $\lambda$ close enough to $HD(\mu)^{-1}$ and $\exp(\chi_\mu)$ respectively, so that
$$
	\frac{\alpha}{\ln \lambda}
		\le \frac{HD(\mu)^{-1}}{\chi_\mu} + \varepsilon_0
		= h_\mu^{-1} + \varepsilon_0.
$$

$6.-$	Fix a big integer $k \ge 1$ and let $V \subset B(z_0, k^{-\alpha})$ be a univalent pull-back of $B(\zeta_0, 3\delta)$ by $f^{\ell(k)}$ and let $x \in V$ be the corresponding preimage of $w_0 \in B(\zeta_0, \delta)$.
	Let $w' \in P^{-k}(x)$ and let $w'' \in B(\zeta_0, \delta)$ a preimage of $w'$ by $P^m$, for some $0 \le m \le m_0$.
	Let $V'$ and $V''$ be the pull-backs of $V$ by $f^{k}$ and $f^{k + m}$ containing $w'$ and $w''$ respectively.
	Since the corresponding pull-backs of $B(z_0, Rk^{-\alpha})$ are univalent, it follows that $\diam(V'') \le \delta$ and therefore $V'' \subset B(\zeta_0, 2\delta)$.  
	Thus $V''$ contains a repelling periodic point of period $k + \ell(k) + m$.
	By Koebe Distortion Theorem there is a universal constant $K_1 > 0$ such that $|(P^{k + \ell(k) + m})'(w'')| \ge K_1 \lambda_{k + \ell(k) + m}$.
	Thus letting $M = \sup_U |P'|$ we have,
\begin{eqnarray*}
	|(P^{k + \ell(k)})'(w')|
		& \ge & M^{-m_0}|(P^{k + \ell(k) + m})'(w'')| \\
		& \ge & K_1 M^{- m_0}\lambda_{k + \ell(k) + m}
		\ge K_1 M^{- m_0}\lambda_{k + \ell(k)}.
\hfill \square
\end{eqnarray*}

\

	Let $\HH \subset \CC$ be the upper half plane and consider a covering map $\psi : \HH \longrightarrow \CC - K(P)$ with deck transformation $z \longrightarrow z + 1$ and such that $P(\psi(z)) = \psi(dz)$.
	In particular $\psi(z + 1) = \psi(z)$, for $z \in \HH$, and for any $r \in \RR$, $\psi$ is injective in $\{ z \in \HH \mid r < {\rm Re}(z) < r + 1 \}$.
	The following are distortion lemmas.

\

\begin{lemma}\label{distortion0}
	There is a constant $D > 1$ such that for all $\tilde{w}$ and $\tilde{w}' \in \HH$ satisfying $|{\rm Im} (\tilde{w})| = |{\rm Im} (\tilde{w}')| \le | \tilde{w} - \tilde{w}'| \le 1/2$,
$$
	|\psi'(\tilde{w})| \le 
		D (|\tilde{w} - \tilde{w}'|/|{\rm Im} (\tilde{w})|)^4 |\psi'(\tilde{w}')|.
$$
\end{lemma}
{\bf Proof.}
	Put $\rho = {\rm Re} (\frac{\tilde{w} + \tilde{w}'}{2})$, $s = |\tilde{w} - \tilde{w}'| \le 1/2$ and $h = {\rm Im}(\tilde{w}) = {\rm Im}(\tilde{w}') \le s$.
	Note that $\psi$ is univalent in the square
$$
	S = \{ z \in \HH \mid |{\rm Re} (z) - \rho| < s 
					\mbox{ and } {\rm Im} (z) < 2s \}
		\subset \{ z \in \HH \mid |{\rm Re} (z) - \rho| < 1/2 \}
$$
	Moreover by hypothesis
$
	\tilde{w}, \tilde{w}' \in S_0 
	= \{ z \in \HH \mid |{\rm Re} (z) - \rho| \le s/2 \mbox{ and } {\rm Im} (z) \le s \}.
$

	Let $\varphi : S \longrightarrow \DD$ be a conformal representation.
	By Schwarz' reflection principle $\varphi$ extends unvalently to $\{ z \in \CC \mid  |{\rm Re}(z) - \rho| < s \mbox{ and } |{\rm Im}(z)| < 2s \}$.
	We normalize $\varphi$ in such a way that $\varphi(\rho) = - i$ and such that the interval $I = (\rho - s, \rho + s) \subset \RR \subset \CC$ is mapped to the semi circle $\{ \exp(i \theta) \mid \pi < \theta < 2 \pi \}$.
	
	Since $\tilde{w}$, $\tilde{w}' \in S_0$ and $I \subset \overline{S}_0$ it follows by Koebe distortion theorem that there is a constant $K > 1$ such that $K^{-1} |\varphi'(\tilde{w})| \le |\varphi'(\tilde{w}')| \le K |\varphi'(\tilde{w})|$ and such that $\widehat{w} = \varphi(\tilde{w})$, $\widehat{w}' = \varphi(\tilde{w}') \in \{ |z| < 1 - K^{-1} s^{-1}h \}$.

	By Theorem 1.3 of \cite{Po} applied to $\psi \circ \varphi^{-1} : \DD \longrightarrow \CC$ there is a universal constant $K_1 > 1$ such that,
$$
	K^{-2} |\psi'(\tilde{w})|/|\psi'(\tilde{w}')|
		\le |(\psi \circ \varphi^{-1})'(\widehat{w})|
			/ |(\psi \circ \varphi^{-1})'(\widehat{w}')|
		\le K_1 (K^{-1}s^{-1}h)^{- 4},
$$
so the lemma follows with $D = K_1 K^6$.
$\hfill \square$

\

\begin{lemma}\label{distortion}
	There is a constant $C_1 = C_1(w_0) > 1$ such that the following property holds.
	Consider an integer $\ell \ge 1$ and $x \in P^{-\ell}(w_0)$.
	Then for every integer $k \ge 1$ and every $w \in P^{-(k + \ell)}(w_0)$ there is $w' \in P^{- k}(x)$ such that,
$$
	|(P^{k + \ell})'(w)| \ge C_1 d^{- 4 \ell} |(P^{k + \ell})'(w')|.
$$ 
\end{lemma}
{\bf Proof.}
	Let $\tilde{w}_0$ and $\tilde{w} \in \HH$ be such that $\psi(\tilde{w}_0) = w_0$ and $\psi(\tilde{w}) = w$.
	We may choose $\tilde{x} \in \HH$ such that $\psi(\tilde{x}) = x$ and $|\tilde{x} - d^k \tilde{w}| \le 1$.
	Put $\tilde{w}' = d^{-k} \tilde{x}$ and $w' = \psi(\tilde{w}') \in P^{-k}(x)$.
	Since ${\rm Im} (\tilde{w}') = d^{-(k + \ell)} {\rm Im}(\tilde{w}_0) = {\rm Im}(\tilde{w})$ the previous lemma implies that $|\psi'(\tilde{w})| \le D {\rm Im}(\tilde{w}_0)^{-4} d^{4 \ell}|\psi'(\tilde{w}')|$.
	
	On the other hand form the equation $P(\psi(z)) = \psi(dz)$ it follows that for every $\widehat{w} \in \HH$ such that $P^{k + \ell}(\psi(\widehat{w})) = w_0$,
$$
	(P^{k + \ell})'(\psi(\widehat{w})) \psi'(\widehat{w}) =
		d^{k + \ell} \psi'(d^{k + \ell} \widehat{w})
		= d^{k + \ell}\psi'(w_0),
$$
which only depends on $w_0$.
	So the lemma follows with $C_1 = D^{-1} {\rm Im}(\tilde{w}_0)^4$.
$\hfill \square$


\
	
	Note that as the invariant measure $\mu$ is chosen with bigger entropy $h_\mu$, the (asymptotic) estimate of Lemma \ref{expansion} is better.
	The topological entropy of $P$ is equal to $\ln d$ (see \cite{Gr} and \cite{MP}), so by the variational principle $h_\mu \le \ln d$; cf. \cite{PU}.
	On the other hand the harmonic measure of $P$ is invariant under $P$ and it has metric entropy equal to $\ln d$, see \cite{Br} and also \cite{FLM} and \cite{Lyu}.
	It follows by Pesin theory that we can choose an invariant measure $\mu$ supported on a hyperbolic set and such that $h_\mu$ is as close to $\ln d$ as wanted, see \cite{PU}.

\

\noindent
{\bf Proof of the Theorem.}
	Put $\lambda_n = Cn^{5 + \varepsilon}$ and let $\varepsilon_0 > 0$ be such that $(4 \ln d) \varepsilon_0 < \varepsilon/3$.
	Let $\mu$ be an invariant probability measure supported on a hyperbolic set and whose metric entropy $h_\mu$ is close enough to $\ln d$ so that $(4 \ln d) h_\mu^{-1} < 4 + \varepsilon/3$.
	Let $w_0 \in \CC - K(P)$, $C_0 = C_0(w_0) > 0$ and $\{ \ell(k) \}_{k \ge 1}$ be given by Lemma \ref{expansion}.
	Moreover let $n \ge 1$ big and $z \in P^{-n}(w_0)$.
	
	Then there is $k$ such that $n \le k + \ell(k) \le n + L + 1$.
	Let $w \in P^{-(k + \ell(k) - n)}(z) \subset P^{-(k + \ell(k))}(w_0)$ and note that $|(P^n)'(z)| \ge M^{-(L + 1)} |(P^{k + \ell(k)})'(w)|$, where $M = \sup |P'|$ taken over all preimages of $w_0$.

	Let $x \in P^{ - \ell(k)}(w_0)$ be given by Lemma \ref{expansion} and let $w' \in P^{- k}(x)$ be given by Lemma \ref{distortion} so that,
$$
	|(P^{k + \ell(k)})'(w)| \ge C_1 d^{- 4 \ell(k)} |(P^{k + \ell(k)})'(w')| \ge C_0C_1 d^{- 4 \ell(k)} \lambda_{k + \ell(k)}.
$$
	By part 1 of Lemma \ref{expansion}, if $k$ is big enough $d^{4 \ell(k)} \le k^{4 \ln d (h_\mu^{-1} + \varepsilon_0)} \le k^{4 + 2 \varepsilon/3}$.
	Thus
$$
	|(P^n)'(z)| \ge M^{-(L + 1)} |(P^{k + \ell(k)})'(w)| 
		\ge M^{-(L + 1)}C_0C_1 \lambda_{k + \ell(k)} k^{- (4 + 2 \varepsilon/3)} \ge C_2 n^{1 + \varepsilon/3},
$$
where $C_2 = M^{-(L + 1)}CC_0C_1$.
	Thus the hypothesis of Lemma \ref{integrable} is satisfied with $\omega_n = C_2n^{1 + \varepsilon/3}$.
$\hfill \square$

\

\subsection*{Acknowledgements}
{\it I'm grateful to F. Przytycki and S. Smirnov for several remarks and comments.}

%
\end{document}